\documentclass[conference,10pt]{IEEEtran}
\usepackage[T1]{fontenc}
\usepackage[latin1]{inputenc}
\usepackage{color}
\usepackage{psfrag}
\usepackage[dvips]{graphicx}
\usepackage{tikz}
\usepackage{pgfplots}
\pgfplotsset{compat=1.12}
\usepackage{enumerate}
\usepackage{rotating}
\usepackage{srcltx}
\usepackage{booktabs}
\usepackage{multicol}
\usepackage{enumitem}
\usepackage{textcomp}
\usetikzlibrary{shapes,arrows}
\usetikzlibrary{arrows,decorations.pathmorphing,backgrounds,positioning,fit,matrix}

\usepackage{xurl}

\usepackage[noadjust]{cite}
\usepackage[cmex10]{amsmath}
\usepackage{amsfonts,bm}
\usepackage{amssymb}
\usepackage{ntheorem}
\usepackage{dsfont}    
\usepackage{mathtools} 
\usepackage{stackengine}

\definecolor{PU_orange}{HTML}{EE7F2D}
\definecolor{PU_darkorange}{HTML}{994400}
\definecolor{PU_lightorange}{HTML}{FFAA66}
\definecolor{PU_black}{HTML}{000000}
\definecolor{PU_darkgray}{HTML}{7F7F83}
\definecolor{PU_lightgray}{HTML}{BDBEC1}

\usepackage{algorithmic}
\usepackage{array}
\usepackage{graphicx}
\usepackage{caption}
\usepackage{subcaption}

\theoremseparator{.}
\newtheorem{corollary}{Corollary}
\newtheorem{prop}{Proposition}
\newtheorem{lem}{Lemma}
\newtheorem{theorem}{Theorem}
\theorembodyfont{\upshape}
\theoremseparator{.}

\newtheorem{remark}{Remark}
\newtheorem*{example*}{Example}

\newcommand{\eu}{\mathrm{e}}

\newcommand{\E}{\mathbb{E}}

\newcommand{\X}{\mathbf{X}}

\newcommand{\Y}{\mathbf{Y}}

\newcommand{\U}{\mathbf{U}}

\newcommand{\sfA}{\mathsf{A}}

\allowdisplaybreaks

\long\def\symbolfootnote[#1]#2{\begingroup%
	\def\thefootnote{\fnsymbol{footnote}}\footnote[#1]{#2}\endgroup} 
\IEEEoverridecommandlockouts    
\title{A General Derivative Identity for the Conditional Expectation with  Focus on the  Exponential Family } 
\author{\IEEEauthorblockN{Alex~Dytso\IEEEauthorrefmark{1},  Martina Cardone\IEEEauthorrefmark{2} }%
\IEEEauthorblockA{\IEEEauthorrefmark{1}New Jersey Institute of Technology, Newark, NJ 07102, USA, Email: alex.dytso@njit.edu \\
\IEEEauthorrefmark{2}University of Minnesota, Minneapolis, MN 55404, USA, Email: mcardone@umn.edu}\\
\vspace{-0.6cm}			   }

\usepackage{geometry}
\begin{document}
\maketitle

\begin{abstract}
Consider a pair of random vectors $(\X,\Y) $  and the conditional expectation operator $\E[\X|\Y={\bf y }]$.  This work studies analytic properties of the conditional expectation by characterizing various derivative identities.  The paper consists of two parts.  
In the first part of the paper,  a general derivative identity for the conditional expectation is derived. Specifically, for the Markov chain $\U \leftrightarrow \X  \leftrightarrow \Y$, a compact expression for  the Jacobian matrix of $\E[\U|\Y = {\bf y}]$ is derived. 
In the second part of the paper, the main identity is specialized to the exponential family.  Moreover, via various choices of the random vector $\U$, the new identity is used to recover and generalize
several known identities and derive some new ones.    
As a first example, a connection between the Jacobian of $ \E[\X|\Y={\bf y }]$ and the conditional variance is established.  As a second example, a recursive expression between higher order conditional expectations is found, which is shown to lead to a generalization of the Tweedy's identity.  
Finally, as a third example, it is shown that the $k$-th order derivative of the conditional expectation is proportional to the $(k+1)$-th order conditional cumulant. 
\end{abstract}  

\section{Introduction}
Consider a pair of random vectors $(\X,\Y) \in \mathbb{R}^{d \times k}$ with a joint distribution $P_{\X \Y}$ and the conditional expectation operator given by
\begin{equation}
\E[\X|\Y={\bf y }]= \int \X \, {\rm d} P_{\X| \Y={\bf y}}. 
\end{equation} 
The conditional expectation operator plays an important role in a variety of fields that require statistical analysis (e.g.,~\cite{banerjee2005optimality,blackwell1947conditional}).
The goal of this work is to study analytical properties of the conditional expectation, i.e., ${\bf y} \mapsto \E[\X|\Y={\bf y }]$. Specifically, the focus is on {\em derivative identities}.  

There exist a number of derivative identities that relate the conditional expectation to other quantities such as the score function, the conditional variance and the conditional cumulants.   
Perhaps the most well-known such identity is the Tweedy's formula~\cite{robbins1956empirical,esposito1968relation}: given that $P_{\Y|\X}$ belongs to the exponential family  with sufficient statistics ${\bf T}({\bf y})$ and base measure $h({\bf y})$, we have that 
\begin{equation}
\left( \boldsymbol{\mathsf{J}}_\mathbf{y}{\bf T}({\bf y}) \right) \E[\X| \Y={\bf y} ] =\nabla_{{\bf{y}}} \log \frac{f_{\Y}({\bf y})}{h({\bf y})},\label{eq:Tweedy_Formula}
\end{equation} 
where  $ \boldsymbol{\mathsf{J}}_\mathbf{y} {\bf T}({\bf y})$ is the Jacobian of the sufficient statistics ${\bf T}({\bf y})$,  and $f_{\Y}(\cdot)$ is the marginal probability density function (pdf) of $\Y$. 
Tweedy's formula is an instrumental tool in statistical signal processing. For example,~\eqref{eq:Tweedy_Formula} implies that the conditional expectation depends on the joint distribution $P_{\X, \Y}$ only through the marginal $P_{\Y}$, which leads to an important class of estimators known as empirical Bayes~\cite{robbins1956empirical}.   For a historic account and impact of Tweedy's formula the interested reader is referred to~\cite{efron2016computer}.   
In information theory, Tweedy's formula is used to connect information and estimator measures and can be used to prove the I-MMSE relationship~\cite{guo2013interplay,palomar2006gradient}.  
The formula is also useful in establishing interesting connections between estimation theory and detection theory~\cite{jaffer1972relations}.

There exist a number of other such identities, which will be surveyed throughout the paper.  However, most of these identities have been derived under some restricted conditions (e.g., $P_{\Y|\X}$ is Gaussian). Moreover, most of these identities are found in an ad-hoc way and there is no unifying approach for characterizing derivative identities for the conditional expectation.  
Recently, in the context of a Gaussian noise model, the authors of~\cite{DytsoISIT2020extended} have provided a unifying approach that is capable of recovering most of the known identities in the literature from a single meta-identity. 
In this work, we extend the results of~\cite{DytsoISIT2020extended} by providing a new general derivative identity that holds under a much more general assumption on the joint distribution $P_{\X, \Y}$. 

\noindent {\bf{Outline and Contributions.}} 
The contributions and the outline of the paper are as follows:
\begin{itemize}[leftmargin=*]
\item In Section~\ref{sec:Main_Idenitty}, Theorem~\ref{thm:MainResult} presents a new identity for the Jacobian of the conditional expectation. Throughout the paper, this general 
identity will be used for systematic proofs of known and novel identities.
\item In Section~\ref{sec:Main_Idenitty_Exponential}, Theorem~\ref{thm:MainIde_Exponential} specializes Theorem~\ref{thm:MainResult} to the case when $P_{\Y|\X}$ is an  exponential family.  Moreover, Theorem~\ref{thm:MainIde_Exponential} is evaluated for the multivariate Gaussian distribution and the Wishart distribution. 
\item Section~\ref{sec:Consequences_Main_Idenity} further focuses on the exponential family and studies consequences of the general derivative identity in Theorem~\ref{thm:MainIde_Exponential}; specifically:
\begin{itemize}
\item In Section~\ref{sec:Variance_idenity}, Proposition~\ref{prop:Variance_Identity}  shows that the Jacobian of $\E[\X|\Y={\bf y }]$ is related to the conditional variance. Moreover, as an application of Proposition~\ref{prop:Variance_Identity}, Section~\ref{sec:Variance_idenity} presents a novel representation of the minimum mean squared error (MMSE).
\item In Section~\ref{sec:Recursive_Identites}, Proposition~\ref{prop:Recursion} focuses on the univariate case and provides a recursive identity between the conditional moments. The recursive identity is used to derive two new results.   First, it is shown that the recursive identity  leads to a new representation of higher-order conditional expectations (i.e., $\E[X^k|Y ]$) in terms of the derivatives of the conditional expectation.
Second, the recursive identity  is used to generalize Tweedy's identity to higher-order conditional expectations. This generalization of the Tweedy's identity maintains  the property that  $\E[X^k|Y]$ depends on the joint distribution only through the marginal of $Y$.
\item In Section~\ref{sec:cumulants},  Theorem~\ref{thm:Cumulants} and Proposition~\ref{prop:Cumulants} establish connections between the derivatives of the conditional cumulant generating function, the derivatives of the conditional cumulants, and the  derivatives of the conditional expectation. For example, it is shown that the $k$-th derivative of the conditional expectation is proportional to the $(k+1)$-th conditional cumulant.
\end{itemize}
\end{itemize}

\noindent {\bf{Notation.}} 
Deterministic scalar quantities are denoted by lowercase letters, scalar random variables are denoted by uppercase letters,  vectors are denoted by bold lowercase letters, random vectors by bold uppercase letters, and matrices by bold uppercase sans serif letters (e.g., $x$, $X$, $\mathbf{x}$, $\mathbf{X}$, $\boldsymbol{\mathsf{X}}$).
$[n_1:n_2]$ is the set of integers from $n_1$ to $n_2 \geq n_1$; $\langle \cdot, \cdot \rangle$ denotes the inner product; for a matrix $\boldsymbol{\mathsf{A}}$, $|\boldsymbol{\mathsf{A}}|$ and ${\rm tr}(\boldsymbol{\mathsf{A}})$ denote the determinant and the trace of $\boldsymbol{\mathsf{A}}$, respectively; $\otimes$ is the Kronecker product; logarithms are in base $\eu$.

For a triplet $(\U,\X,\Y) \in \mathbb{R}^{m \times n \times k}$ we define the \emph{conditional variance matrix} and the \emph{conditional cross-covariance matrix} as follows,
\begin{align}
 \boldsymbol{\mathsf{Var}}(\X | \Y) &= \E \left[  \X   \X^\mathsf{T}   |\Y\right ]    -  \E \left[  \X   |\Y \right ]  \E \left[  \X^\mathsf{T} | \Y\right], \label{eq:CondVar}\\
\boldsymbol{\mathsf{Cov}}(\X, \U | \Y)  &= \E \left[  \X   \U^\mathsf{T}   |\Y\right ]    -  \E \left[  \X   |\Y \right ]  \E \left[  \U^\mathsf{T} | \Y\right] . \label{eq:CondCovar}
\end{align} 
The \emph{gradient} of a  function $f: \mathbb{R}^n \to \mathbb{R}$  is denoted by 
 \begin{equation}
 \nabla_{\bf x} f(\mathbf{x})  = \begin{bmatrix}     \frac{ \partial  f({\bf x}) }{ \partial  x_1}  &\frac{  \partial   f({\bf x}) }{ \partial  x_2}   & \ldots &  \frac{ \partial   f({\bf x})}{ \partial  x_n}   \end{bmatrix}^\mathsf{T} \in \mathbb{R}^n .
 \end{equation}
 The \emph{Jacobian matrix}  of a function $ \mathbf{f}: \mathbb{R}^n \to \mathbb{R}^m$ is denoted by $\boldsymbol{\mathsf{J}}_\mathbf{x} \mathbf{f}(\mathbf{x}) \in \mathbb{R}^{n \times m}$ and defined as
 \begin{equation}
 \label{eq:Def_Jacobian}
 \boldsymbol{\mathsf{J}}_\mathbf{x} \mathbf{f}(\mathbf{x}) =  
 \begin{bmatrix}
 \nabla_{\bf x} {f}_1(\mathbf{x}) & \nabla_{\bf x} {f}_2(\mathbf{x}) & \ldots & \nabla_{\bf x} {f}_m(\mathbf{x})
 \end{bmatrix} .
 \end{equation}
 %
  The {\em vectorization} of  a matrix $\boldsymbol{\mathsf{A}}  \in  \mathbb{R}^{n \times m}$ is given by
  \begin{align}
  \label{eq:Vec}
  & \boldsymbol{vec}(\boldsymbol{\mathsf{A}}) = \begin{bmatrix} \mathbf{a}_1^\mathsf{T} &  \mathbf{a}_2^\mathsf{T} & \ldots & \mathbf{a}_m^\mathsf{T} \end{bmatrix}^\mathsf{T} \in \mathbb{R}^{m n \times 1},
     \end{align}
where, for $i \in [1:m]$, $\mathbf{a}_i \in \mathbb{R}^{n \times 1}$ is the $i$-th column of~$\boldsymbol{\mathsf{A}}$.
     
     For a symmetric matrix  $\boldsymbol{\mathsf{A}}  \in  \mathbb{R}^{n \times n}$, the {\em half-vectorization operator} is defined as the vectorization of the lower triangle part of $\boldsymbol{\mathsf{A}} $, i.e., $ \boldsymbol{vech}(\boldsymbol{\mathsf{A}}) \in \mathbb{R}^{\frac{n(n+1)}{2} \times 1}$ is defined as
\begin{align}
\label{eq:HalfVec}
& \boldsymbol{vech}(\boldsymbol{\mathsf{A}}) = \begin{bmatrix} \mathbf{a}_1^\mathsf{T}[1:n] &  \mathbf{a}_2^\mathsf{T}[2:n]  & \ldots  &  \mathbf{a}_n^\mathsf{T}[n:n] \end{bmatrix}^\mathsf{T} ,
\end{align}     
where we use the (Matlab-inspired) notation $\mathbf{a}_i[i:n] \in \mathbb{R}^{(n-i+1)\times 1}$ to denote the $i$-th column of $\boldsymbol{\mathsf{A}}$ where only the elements from the $i$-th row to the $n$-th row are retained.
%
%
   Finally, the duplication matrix and the elimination matrix are denoted by $\boldsymbol{\mathsf{D}}_n \in \mathbb{R}^{n^2 \times \frac{n(n+1)}{2}}$ and $\boldsymbol{\mathsf{L}}_n \in \mathbb{R}^{ \frac{n(n+1)}{2} \times n^2}$, respectively, and satisfy
   \begin{equation}
   \label{eq:ElDupMat}
  \boldsymbol{\mathsf{D}}_n \boldsymbol{vech}(\boldsymbol{\mathsf{A}})= \boldsymbol{vec}(\boldsymbol{\mathsf{A}}), \quad \boldsymbol{\mathsf{L}}_n   \boldsymbol{vec}(\boldsymbol{\mathsf{A}})=\boldsymbol{vech}(\boldsymbol{\mathsf{A}}).
   \end{equation}
   The space of symmetric positive definite matrices of dimension $n$ is denoted by $\mathbb{S}^n_{+}$. 

\section{Main Derivative Identity}
\label{sec:Main_Idenitty} 

In order to present the main derivative identity and its proof, we need the notions of information density, score function and conditional score function. The {\em information density} for the joint distribution $P_{\X,\Y}$  supported on $\mathcal{X} \times \mathcal{Y}$ where $\mathcal{X} \subseteq \mathbb{R}^d$ and $\mathcal{Y} \subseteq \mathbb{R}^k$  is defined as
\begin{equation}
\label{eq:InfoDens}
\iota_{P_{\X,\Y}}( {\bf x};{\bf y})=\log \frac{{\rm d} P_{\X,\Y}}{{\rm d} ( P_{\X}  \cdot P_{\Y}) } ({\bf x}, {\bf y} ),\,    {\bf x} \in \mathcal{X}, \,  {\bf y}  \in\mathcal{Y},
\end{equation}
where $\frac{{\rm d} P_{\X,\Y}}{{\rm d} ( P_{\X}  \cdot P_{\Y}) } ({\bf x}, {\bf y} )= \frac{{\rm d} P_{\Y| \X ={\bf x}}}{{\rm d}  P_{\Y} } ({\bf x}, {\bf y} ) $ is the Radon-Nikodym derivative with the understanding that  if $P_{\Y| \X ={\bf x}}$ is not absolutely continuous with respect to $P_{\Y}$ we define $\iota_{P_{\X,\Y}}( {\bf x};{\bf y})=\infty$.
For a distribution $P_\Y$ with  a pdf $f_{\Y}(\cdot)$, the {\em score function} is defined as
\begin{equation}
\rho_{\Y}({\bf y})=\nabla_{\bf y} \log f_{\Y}({\bf y})=\frac{\nabla_{\bf y}  f_{\Y}({\bf y})}{ f_{\Y}({\bf y})},\, {\bf y} \in \mathcal{Y}. 
\end{equation}
For the conditional distribution $P_{\Y|\X}$ with a pdf  $f_{\Y|\X}(\cdot | \cdot)$, we define the {\em conditional score function} as
\begin{equation}
\label{eq:ConScore}
\rho_{\Y|\X}({\bf y}|{\bf x})= \nabla_{\bf y}  \log f_{\Y|\X}({\bf y}| {\bf x} ),  \,    {\bf x} \in \mathcal{X}, \,  {\bf y}  \in\mathcal{Y}. 
\end{equation}
The following elementary properties of the above three quantities are proved for completness in Appendix~\ref{app:ProofLemmaRelations} and  will be useful in the proof of the main identity.
\begin{lem} \label{lem:Prop}
The information density, score function, and conditional score function satisfy the following properties:
\begin{itemize}[leftmargin=*]
\item (Conditional Expectation and Information Density):  Suppose that $ \U \leftrightarrow \X \leftrightarrow \Y$  forms a Markov chain, in that order. Then,  for ${\bf y} \in \mathcal{Y}$
\begin{equation}
\E[\U|\Y={\bf y} ]=  \E \left [\U \eu^{\iota_{P_{\X,\Y}}( \X ;{\bf y})} \right ]. \label{eq:CE_info_density}
\end{equation}
\item  (Gradient of the Information Density): Suppose  that the distributions $P_\Y$ and $P_{\Y|\X}$ have pdfs $f_{\Y}(\cdot)$ and $f_{\Y|\X}(\cdot | \cdot)$, respectively. Then,  for $({\bf x}, {\bf y}) \in   \mathcal{X} \times \mathcal{Y}$
\begin{equation}
\nabla_{\bf y}\iota_{P_{\X,\Y}}( {\bf x} ;{\bf y})=  \rho_{\Y|\X}({\bf y}|{\bf x})-\rho_{\Y}({\bf y}).   \label{eq:Gradient_Info_density}
\end{equation}
\item  (Score Function vs. Conditional Score Function): Suppose  that $\E \left[  \| \nabla_{\bf y} f_{\Y|\X}({\bf y}|  \X ) \| \right]<\infty$ for all ${\bf y} \in \mathcal{Y}$ and that $({\bf x},{\bf y})  \mapsto f_{\Y|\X}({\bf y}|{\bf x})$ is absolutely continuous in ${\bf y} $ for every ${\bf x}$. Then, for ${\bf y} \in \mathcal{Y}$
\begin{equation}
\rho_{\Y}({\bf y})= \E [\rho_{\Y|\X}(\Y|\X)|\Y={\bf y} ].\label{eq:Scond_average}
\end{equation}
\end{itemize} 
\end{lem} 
By leveraging the properties in Lemma~\ref{lem:Prop}, we can now prove the main derivative identity, which is provided in the next theorem.
\begin{theorem} \label{thm:MainResult} Suppose that the random vectors  $\U \in \mathbb{R}^m, \X \in \mathcal{X} \subseteq  \mathbb{R}^d$ and $\Y \in \mathcal{Y} \subseteq \mathbb{R}^k$ satisfy the following conditions:
\begin{align*}
{\rm{A1:}}& \, \U \leftrightarrow \X \leftrightarrow \Y \text{ forms a Markov chain, in that order;}\\
{\rm{A2:}}& \, \E \left[ \left \| \U \right \|   \left \| \nabla_{ \bf y} \iota_{P_{\X,\Y}}( \X ; \Y) \right \| | \Y={\bf y}  \right ]<\infty, \, {\bf y} \in \mathcal{Y}  \subseteq \mathbb{R}^k;\\
{\rm{A3:}} & \, \text{$\E \left[  \| \nabla_{\bf y} f_{\Y|\X}({\bf y}|  \X ) \| \right]<\infty$ for all ${\bf y} \in \mathcal{Y}$};\\
{\rm{A4:}} & \,  \text{$({\bf x},{\bf u},{\bf y})  \mapsto {\bf u} \, \iota_{P_{\X,\Y}}( {\bf x} ; {\bf y})$ and $({\bf x},{\bf y})  \mapsto f_{\Y|\X}({\bf y}|{\bf x})$ }  \\
& \text{ are absolutely continuous in ${\bf y} $ for every $({\bf x},{\bf u})$}. 
\end{align*}
Then, for ${\bf y} \in \mathcal{Y}$
\begin{equation}
\label{eq:IdentGen}
\boldsymbol{\mathsf{J}}_\mathbf{y} \E[\U|\Y={\bf y}]= \boldsymbol{\mathsf{Cov}} \left( \hspace{0.01cm}  \rho_{\Y|\X}(\Y|\X), \U \, | \, \Y={\bf y} \right).
\end{equation}
\end{theorem} 
\begin{IEEEproof}
For $i \in [1:m] $, we have that
\begin{align}
&\nabla_{ \bf y} \E[ U_i|\Y={\bf y}] \notag\\
&\stackrel{\rm{(a)}}{=}  \nabla_{ \bf y}  \E \left[ U_i \eu^{\iota_{P_{\X,\Y}}( \X ;{\bf y})} \right ] \notag \\
&\stackrel{\rm{(b)}}{=}   \E \left[ U_i  \nabla_{ \bf y} \iota_{P_{\X,\Y}}( \X ;{\bf y})   \eu^{\iota_{P_{\X,\Y}}( \X ;{\bf y})} \right ] \notag \\
&\stackrel{\rm{(c)}}{=}   \E \left[ U_i  \nabla_{ \bf y} \iota_{P_{\X,\Y}}( \X ;{\bf y})    |\Y={\bf y} \right]\notag \\
&\stackrel{\rm{(d)}}{=}   \E[ U_i  ( \rho_{\Y|\X}(\Y|\X)- \rho_{\Y}(\Y))    |\Y={\bf y}] \notag \\
&=   \E[ U_i   \rho_{\Y|\X}(\Y|\X)  |\Y={\bf y}] -   \E[ U_i     |\Y={\bf y}] \rho_{\Y}({\bf y}) \notag \\
&\stackrel{\rm{(e)}}{=}  \boldsymbol{\mathsf{Cov}} \left(  \rho_{\Y|\X}(\Y|\X), U_i | \Y={\bf y}  \right), \label{eq:Coordinate_i_result}
\end{align}
where the labeled equalities follow from: $\rm{(a)}$ using~\eqref{eq:CE_info_density}, which holds under the assumption ${\rm{A1}}$;   ${\rm{(b)}}$   interchanging the gradient and expectation that is permissible by using  the Leibniz integral rule, which requires verifying the conditions in ${\rm{A2}}$ and ${\rm{A4}}$; ${\rm{(c)}}$ using~\eqref{eq:CE_info_density}, which holds under ${\rm{A1}}$;  
${\rm{(d)}}$ using~\eqref{eq:Gradient_Info_density}; and  ${\rm{(e)}}$ using~\eqref{eq:Scond_average}, which holds under the assumptions ${\rm{A3}}$ and ${\rm{A4}}$, and the definition of  conditional covariance in~\eqref{eq:CondCovar}. 

The proof of Theorem~\ref{thm:MainResult} is concluded by using~\eqref{eq:Coordinate_i_result} together with the definition of  the Jacobian in~\eqref{eq:Def_Jacobian}. 
\end{IEEEproof} 

The identity in~\eqref{eq:IdentGen} has a number of interesting consequences. In the remaining of the paper, we explore these consequences in the context of an exponential family.

\section{Exponential Family and the Main Identity}
\label{sec:Main_Idenitty_Exponential}
The class of probability models $\mathcal{P}=\{ P_{ \Y|\X={\bf x}},  {\bf x} \in \mathcal{X} \subseteq \mathbb{R}^d \}$ supported on $\mathcal{Y} \subseteq \mathbb{R}^k$ is an \emph{exponential family} if the pdf of it can be written as
\begin{equation}
\label{eq:ExpFam}
f_{\Y|\X}({\bf y}| {\bf x}) 
=h({\bf{y}}) \eu^{ \langle {\bf{x}}, {\bf T}({\bf y}) \rangle -\phi({\bf{x}}) }, \, {\bf y} \in \mathcal{Y}, \,  
{\bf{x}} \in  \mathcal{X}, 
\end{equation}
where   
${\bf T}:\mathcal{Y} \to \mathbb{R}^d $ is \emph{the sufficient statistic function}; 
$\phi: \mathcal{X} \to \mathbb{R}$ is \emph{the log-partition function}; 
and $h: \mathcal{Y} \to [0,\infty)$ is \emph{the base measure}.   
In this work, we assume that $\mathcal{Y}$ is an open set and that the base measure $h$ is absolutely continuous with respect to the Lebesgue measure.  
In other words, we restrict our focus to continuous distributions belonging to the exponential family (e.g., normal, Wishart). Furthermore, in the remaining, we assume that the sufficient statistics ${\bf T}({\bf y})$ is an analytic function.

Note that there is a number of known derivative identities for the mean and variance of the exponential distribution. For example,   the mean  and the  variance can be related to the log-partition function $\phi(\cdot)$ via the following relationships: 
\begin{equation}
\E[\Y|\X={\bf x}] \!=\!\nabla \phi({\bf x}) ,   \quad  \boldsymbol{\mathsf{Var}}(\Y| \X={\bf x} )\!=\! {\bf H}^2   \phi({\bf x}),  \label{eq:Some_Known_Derivative_Identities}
\end{equation}
where ${\bf H}^2$ is the Hessian matrix of $\phi(\cdot)$.   
Note that the  identities in \eqref{eq:Some_Known_Derivative_Identities} do   not assume any prior distribution on $\X$.  In this work, however,  we are interested in a family of derivative identities that are fundamentally different from the above identities. Specifically, we are interested in studying derivative identities for quantities such  as $\E[\X|\Y={\bf y}]$ in which there is a prior distribution on $\X$.  

The next theorem specializes Theorem~\ref{thm:MainResult} to the aforementioned exponential family.

%
%
\begin{theorem}  
\label{thm:MainIde_Exponential}
Suppose that random vectors  $\U \in \mathbb{R}^m, \X \in \mathcal{X} \subseteq  \mathbb{R}^d$ and $\Y \in \mathcal{Y} \subseteq \mathbb{R}^k$ are such that $f_{\Y|\X}(\cdot| \cdot)$ is as in~\eqref{eq:ExpFam}, and satisfy the following conditions,
\begin{align*}
{\rm{A1:}}&\U \leftrightarrow \X \leftrightarrow \Y \text{ forms a Markov chain, in that order;}\\
{\rm{A2:}}& \E \left [\left \| \U \right \|  | \Y={\bf y}\right ] < \infty, \, {\bf y} \in \mathcal{Y}  \subseteq \mathbb{R}^k;\\
{\rm{A3:}}&\E \left[ \left \| \U \right \|   \left \|  \X \right \| | \Y={\bf y}  \right ]<\infty, \, {\bf y} \in \mathcal{Y}  \subseteq \mathbb{R}^k;\\
{\rm{A4:}} & \,  \text{$({\bf x},{\bf u},{\bf y})  \mapsto {\bf u} \, \iota_{P_{\X,\Y}}( {\bf x} ; {\bf y})$ and $({\bf x},{\bf y})  \mapsto f_{\Y|\X}({\bf y}|{\bf x})$ }  \\
& \text{ are absolutely continuous in ${\bf y} $ for every $({\bf x},{\bf u})$}; \\
{\rm{A5:}}&\E \left[   \left \|   \X \right \| | \Y={\bf y}  \right ]<\infty,  \, {\bf y} \in \mathcal{Y}  \subseteq \mathbb{R}^k. 
\end{align*} 
Then,  for $ {\bf y} \in \mathcal{Y}  \subseteq \mathbb{R}^k$, we have that
\begin{equation}
\label{eq:IdentExp}
\boldsymbol{\mathsf{J}}_\mathbf{y} \E[\U|\Y={\bf y}]=  \left( \boldsymbol{\mathsf{J}}_\mathbf{y}  {\bf T}({\bf y})  \right) \boldsymbol{\mathsf{Cov}}(\X,  \U  | \Y={\bf y})  .
\end{equation}
\end{theorem} 
\begin{IEEEproof}
In Appendix~\ref{app:AssumExpFamily}, we show how the three conditions in Theorem~\ref{thm:MainResult} specialize to the case of the exponential family in~\eqref{eq:ExpFam}.
Now, we are left to show the identity in~\eqref{eq:IdentExp}. From~\eqref{eq:IdentGen}, we obtain
\begin{align*}
\boldsymbol{\mathsf{J}}_\mathbf{y} \E[\U|\Y={\bf y}]&= \boldsymbol{\mathsf{Cov}} \left(  \rho_{\Y|\X}(\Y|\X), \U | \Y={\bf y} \right)
\\ & \stackrel{{\rm{(a)}}}{=}\boldsymbol{\mathsf{Cov}} \left(   \nabla_{\bf y} \log (h({\bf y}))  , \U | \Y={\bf y} \right)
\\& \quad + \boldsymbol{\mathsf{Cov}} \left(   \left( \boldsymbol{\mathsf{J}}_\mathbf{y}  {\bf T}({\bf y})  \right) \X , \U | \Y={\bf y} \right)
\\ & \stackrel{{\rm{(b)}}}{=} \left( \boldsymbol{\mathsf{J}}_\mathbf{y}  {\bf T}({\bf y})  \right) \boldsymbol{\mathsf{Cov}} \left(    \X , \U | \Y={\bf y} \right),
\end{align*}
where $\rm{(a)}$ follows from the fact that $\boldsymbol{\mathsf{Cov}}(\mathbf{A}+\mathbf{B},\mathbf{C}) = \boldsymbol{\mathsf{Cov}}(\mathbf{A},\mathbf{C}) + \boldsymbol{\mathsf{Cov}}(\mathbf{B},\mathbf{C})$ and the fact that, for the exponential family in~\eqref{eq:ExpFam}, the conditional score function is given by 
\begin{align}
\label{eq:CondScorExp}
 \rho_{\Y|\X}(\Y|\X)&= \nabla_{\bf y} \log (h({\bf y})) +  \nabla_{\bf y} \langle {\bf{x}}, {\bf T}({\bf y}) \rangle \notag \\
& = \nabla_{\bf y} \log (h({\bf y})) + \left( \boldsymbol{\mathsf{J}}_\mathbf{y}  {\bf T}({\bf y})  \right)   \bf{x};
 \end{align}
and $\rm{(b)}$ holds since $\boldsymbol{\mathsf{Cov}} \left(   \nabla_{\bf y} \log (h({\bf y}))  , \U | \Y={\bf y} \right) =\mathsf{\bf 0}$.
This concludes the proof of Theorem~\ref{thm:MainIde_Exponential}.
\end{IEEEproof} 

\begin{remark}
The assumptions ${\rm{A2}}$, ${\rm{A3}}$ and ${\rm{A5}}$ in Theorem~\ref{thm:MainIde_Exponential} are rather mild. For example, these assumptions hold if both $\X$ and $\U$ are integrable, i.e., $\E[\| \X\|^2 ]<\infty, \E[\| \U\|^2 ]<\infty$. In the remaining of the paper, we only consider priors on $(\U,\X)$ that satisfy the regularity conditions ${\rm{A1}}-{\rm{A5}}$.
\end{remark}

We  conclude this section by evaluating the result in Theorem~\ref{thm:MainIde_Exponential} for two continuous distributions that belong to the exponential family.

\begin{example*} \emph{Multivariate Normal with (Un)Known Mean and (Un)Known Covariance.}
Consider the case where $\Y | \X={\bf x} \sim \mathcal{N}( \bm{m};  \boldsymbol{\mathsf{\Sigma}})$  
where  $\mathcal{N}( \bm{m};  \boldsymbol{\mathsf{\Sigma}})$  is the multivariate normal distribution with mean $\bm{m}$ and covariance matrix $\boldsymbol{\mathsf{\Sigma}}$. The mapping to the exponential family in \eqref{eq:ExpFam} is done as follows,
\begin{subequations}
\begin{align}
{\bf x}&=\begin{bmatrix}  \left [ \boldsymbol{\mathsf{\Sigma}}^{-1} \bm{m} \right]^\mathsf{T} & \left [ \boldsymbol{vec} \left(-\frac{1}{2} \boldsymbol{\mathsf{\Sigma}}^{-1} \right) \right ]^\mathsf{T}\end{bmatrix}^\mathsf{T}, \\
h({\bf y})&=(2\pi)^{-k/2},\\
 {\bf T}({\bf y}) &= \begin{bmatrix}{\bf y}^\mathsf{T} &  \left[\boldsymbol{vec} \left( {\bf y} {\bf y}^\mathsf{T} \right)\right ]^\mathsf{T} \end{bmatrix}^\mathsf{T}, \\
 \phi({\bf x})&=\frac{1}{2}  \left(  \bm{m}^\mathsf{T} \boldsymbol{\mathsf{\Sigma}}^{-1}    \bm{m} +\log(|\boldsymbol{\mathsf{\Sigma}} |) \right).
\end{align}
\end{subequations}
The Jacobian of ${\bf T}({\bf y}) $ is then given by
\begin{equation*}
\boldsymbol{\mathsf{J}}_\mathbf{y} {\bf T}({\bf y}) = \begin{bmatrix} \boldsymbol{\mathsf{I}}_k  & \mathbf{y}^\mathsf{T}   \otimes \boldsymbol{\mathsf{I}}_k  +  \boldsymbol{\mathsf{I}}_k  \otimes \mathbf{y}^\mathsf{T} \end{bmatrix}. 
\end{equation*}
Now, suppose that $\X=\begin{bmatrix}  \left [ \boldsymbol{\mathsf{\Sigma}}^{-1} \bm{M} \right]^\mathsf{T} & \left [ \boldsymbol{vec} \left(-\frac{1}{2} \boldsymbol{\mathsf{\Sigma}}^{-1} \right) \right ]^\mathsf{T}\end{bmatrix}^\mathsf{T}$ where the joint $ P_{ \bm{M}, \boldsymbol{\mathsf{\Sigma}}}$ satisfies  $\rm{A1}$-$\rm{A5}$ of Theorem~\ref{thm:MainIde_Exponential}. Then,  for ${\bf y} \in \mathbb{R}^k$, from~\eqref{eq:IdentExp} we~obtain
\begin{align}
\label{eq:IdentityMultGaussi}
&\boldsymbol{\mathsf{J}}_\mathbf{y} \E[\U|\Y={\bf y}] =    \begin{bmatrix}\boldsymbol{\mathsf{I}}_k,  \, \, \mathbf{y}^\mathsf{T}   \otimes \boldsymbol{\mathsf{I}}_k  +  \boldsymbol{\mathsf{I}}_k  \otimes \mathbf{y}^\mathsf{T} \end{bmatrix} \notag
\\&\qquad \qquad \quad \cdot \boldsymbol{\mathsf{Cov}}\left( \begin{bmatrix} \boldsymbol{\mathsf{\Sigma}}^{-1} \bm{M} \\ 
\boldsymbol{vec} \left(-\frac{1}{2}\boldsymbol{\mathsf{\Sigma}}^{-1} \right)
 \end{bmatrix}   ,  \U   \Big | \Y={\bf y}    \right) .
\end{align}
The above can now be specialized in two ways. 
First, specializing it to the case where $\boldsymbol{\mathsf{\Sigma}}$ is known  and $\bm{M}$ is unknown (i.e., $P_{  \boldsymbol{\mathsf{\Sigma}}}$ is a point mass), we arrive at
\begin{align*}
&\boldsymbol{\mathsf{J}}_\mathbf{y} \E[\U|\Y={\bf y}] =    \boldsymbol{\mathsf{\Sigma}}^{-1}  \boldsymbol{\mathsf{Cov}}\left(    \bm{M} , \U \Big | \Y={\bf y} \right),
\end{align*}
which recovers the result in \cite{DytsoISIT2020extended}.   
Second, specializing~\eqref{eq:IdentityMultGaussi} to the case when $\bm{M}$ is known and equal to zero but  $\boldsymbol{\mathsf{\Sigma}}$ is  unknown, we arrive at
\begin{align*}
&\boldsymbol{\mathsf{J}}_\mathbf{y} \E[\U|\Y={\bf y}] 
\\&= - \frac{1}{2} \left[  
\mathbf{y}^\mathsf{T}   \otimes \boldsymbol{\mathsf{I}}_k  +  \boldsymbol{\mathsf{I}}_k  \otimes \mathbf{y}^\mathsf{T}
 \right]  \boldsymbol{\mathsf{Cov}}\left(  
\boldsymbol{vec} \left( \boldsymbol{\mathsf{\Sigma}}^{-1} \right), \U
  \Big | \Y={\bf y} \right)      .
\end{align*}
\end{example*}


\begin{example*} 
{\em Wishart distribution.}
The Wishart distribution with $n$ degrees of freedom and a parameter matrix $ \boldsymbol{\mathsf{V}}  \in \mathbb{S}^p_{+}$, where $n \ge p$, is given by
\begin{equation}
f_{\boldsymbol{\mathsf{\sfA}}|\boldsymbol{\mathsf{V}},N}(\boldsymbol{\mathsf{\sfA}}|\boldsymbol{\mathsf{V}},n)=    \frac{| \boldsymbol{\mathsf{\sfA}} |^{\frac{n-p-1}{2}} \eu^{-\frac{{\rm tr}(\boldsymbol{\mathsf{V}}^{-1} \boldsymbol{\mathsf{\sfA}} )}{2}} }{2^{\frac{np}{2}} \Gamma_p \left(  \frac{n}{2}\right) | \boldsymbol{\mathsf{V}} |^{\frac{n}{2}}  }, \, \boldsymbol{\mathsf{\sfA}}  \in \mathbb{S}^p_{+} ,
\end{equation}
where $ \Gamma_p( \cdot)$ is the multivariate gamma function. 
The pdf of the Wishart distribution can be written in the exponential form in \eqref{eq:ExpFam} by using the following mappings\footnote{Since $\boldsymbol{\mathsf{A}}$ belongs to the set of symmetric matrices, it is more convenient to take the Jacobian in~\eqref{eq:IdentExp} with respect to $\boldsymbol{vech}(\boldsymbol{\mathsf{A}})$ instead of $\boldsymbol{vec}(\boldsymbol{\mathsf{A}})$.},
%
\begin{subequations}
\begin{align}
{\bf y}&=\boldsymbol{vech}(\boldsymbol{\mathsf{A}}), \label{eq:yWish}\\
{\bf x}&= \begin{bmatrix}-\frac{1}{2}  \left [ \boldsymbol{vec}(\boldsymbol{\mathsf{V}}^{-1}) \right ]^\mathsf{T} & \frac{n-p-1}{2} \end{bmatrix}^\mathsf{T},  \\
h( {\bf y})&=1, \\
{\bf T}({\bf y})&=  \begin{bmatrix} \left [ \boldsymbol{vec}(\boldsymbol{\mathsf{A}}) \right ]^\mathsf{T} & \log( |\boldsymbol{\mathsf{A}}|) \end{bmatrix}^\mathsf{T},\\
\phi({\bf x})&=\frac{n}{2} \log (| \boldsymbol{\mathsf{V}}|)+\log \Gamma_p \left(  \frac{n}{2}\right) + \frac{np}{2} \log(2).
\end{align} 
\end{subequations}
We also note that the Jacobian of ${\bf T}({\bf y}) $ is given by 
\begin{align*}
\boldsymbol{\mathsf{J}}_\mathbf{y} {\bf T}({\bf y}) = \begin{bmatrix} \boldsymbol{\mathsf{D}}_p^\mathsf{T} & 
\boldsymbol{\mathsf{D}}_p^\mathsf{T} \boldsymbol{\mathsf{D}}_p  \mathbf{y} \end{bmatrix},
\end{align*} 
where we have used the facts that
\begin{align*}
\boldsymbol{\mathsf{J}}_\mathbf{y} \boldsymbol{vec}(\boldsymbol{\mathsf{A}})&\stackrel{\eqref{eq:ElDupMat}}{=}\boldsymbol{\mathsf{J}}_\mathbf{y} \left( \boldsymbol{\mathsf{D}}_n \boldsymbol{vech}(\boldsymbol{\mathsf{A}}) \right ) \stackrel{\eqref{eq:yWish}}{=} \boldsymbol{\mathsf{J}}_\mathbf{y}  \left( \boldsymbol{\mathsf{D}}_p  \mathbf{y} \right )= \boldsymbol{\mathsf{D}}_p^\mathsf{T} ,\\
\nabla_\mathbf{y} \log( |\boldsymbol{\mathsf{A}}|) &= 
\frac{\left. \boldsymbol{\mathsf{D}}_p^\mathsf{T} \boldsymbol{\mathsf{J}}_{\boldsymbol{vec}(\boldsymbol{\mathsf{A}})} (|\boldsymbol{\mathsf{A}}|) \right |_{\boldsymbol{vec}(\boldsymbol{\mathsf{A}}) =  \boldsymbol{\mathsf{D}}_p  \mathbf{y}}}{|\boldsymbol{\mathsf{A}}|} 
\stackrel{\eqref{eq:yWish}}{=}  \boldsymbol{\mathsf{D}}_p^\mathsf{T} \boldsymbol{\mathsf{D}}_p  \mathbf{y}.
\end{align*} 
Next, by assuming a prior distribution on the parameter 
${\bf X}= \begin{bmatrix}-\frac{1}{2}  \left [ \boldsymbol{vec}(\boldsymbol{\mathsf{V}}^{-1}) \right ]^\mathsf{T} & \frac{N-p-1}{2} \end{bmatrix}^\mathsf{T}$, where it is assumed that $(\boldsymbol{\mathsf{V}}, N) \sim P_{\boldsymbol{\mathsf{V}},N}$, the identity in~\eqref{eq:IdentExp} reduces to 
\begin{align}
\label{eq:MainIdWish}
&\boldsymbol{\mathsf{J}}_\mathbf{y} \E[\U|\Y={\bf y}] =\begin{bmatrix} \boldsymbol{\mathsf{D}}_p^\mathsf{T} & 
\boldsymbol{\mathsf{D}}_p^\mathsf{T} \boldsymbol{\mathsf{D}}_p  \mathbf{y} \end{bmatrix} \notag 
\\&\qquad \qquad \quad \cdot  \boldsymbol{\mathsf{Cov}}\left(     \left[   \begin{array}{l} -\frac{1}{2}  \boldsymbol{vec}(\boldsymbol{\mathsf{V}}^{-1}) \\ 
\frac{N-p-1}{2} 
 \end{array}  \right] ,  \U   \Big | \Y={\bf y}    \right). 
\end{align}
We next specialize the above result to the case $p=1$, i.e., the gamma distribution. In other words, we consider
\begin{align*}
f_{Y|A,B}( y| a,b) = \frac{ y^{a-1} \eu^{- y b}}{b^{-a} \Gamma(a)}, \quad y,a,b >0.
\end{align*}
For this case, the result in~\eqref{eq:MainIdWish} reduces to
\begin{align*}
\frac{{\rm d}}{{\rm{d}} y} \E[\U |Y=y] = \begin{bmatrix} 1 & y\end{bmatrix}
\boldsymbol{\mathsf{Cov}}\left(
\begin{bmatrix}
-B
\\
A
\end{bmatrix},
\U
 \Big | Y=y \right ).
\end{align*}

\end{example*}

\section {Consequences of the Main Identity for the Exponential Family} 
\label{sec:Consequences_Main_Idenity}
In this section, we show that several well-known identities for the exponential family defined in~\eqref{eq:ExpFam} can be derived systematically from the Jacobian identity in Theorem~\ref{thm:MainIde_Exponential}. 
Moreover, we use this new identity to derive several generalizations of previously known identities and discover some new identities. 
Specifically, we will evaluate Theorem~\ref{thm:MainIde_Exponential} with three different choices of~$\U$.

\subsection{Variance Identity}
\label{sec:Variance_idenity}
By setting  $\U=\X$ in Theorem~\ref{thm:MainIde_Exponential} we arrive at the following result. 
\begin{prop}\label{prop:Variance_Identity} 
For the exponential family defined in~\eqref{eq:ExpFam}, we have that
\begin{equation}
\boldsymbol{\mathsf{J}}_\mathbf{y} \E[\X|\Y={\bf y}]= \left( \boldsymbol{\mathsf{J}}_\mathbf{y}  {\bf T}({\bf y})  \right)   \boldsymbol{\mathsf{Var}}(  \X  | \Y={\bf y}).  \label{eq:Variance_Identity} 
\end{equation}
\end{prop} 
The identity in~\eqref{eq:Variance_Identity} was previously demonstrated for the case of a Gaussian noise channel with known variance. 
Specifically, in~\cite{hatsell1971some} this identity was proven for the vector Gaussian noise with an identity covariance matrix, and then generalized to an arbitrary covariance matrix~\cite{palomar2006gradient}.  
We refer the interested reader to~\cite{DytsoISIT2020extended} for an account of the impact of this identity in the case of Gaussian noise.

Recall that the MMSE matrix  of estimating  $\X$ from $\Y$ is given by 
\begin{equation*}
\boldsymbol{\mathsf{MMSE}}(\X| \Y)=   \E \left[  (\X       -  \E \left[  \X   |\Y \right ]) (\X       -  \E \left[  \X   |\Y \right ])^\mathsf{T}   \right].
\end{equation*} 
As an application of Proposition~\ref{prop:Variance_Identity}, we have the following new representation of the MMSE matrix. 
\begin{corollary} 
Suppose that $\left( \boldsymbol{\mathsf{J}}_\mathbf{y}  {\bf T}({\bf y})  \right) $ is $P_{\Y}$-almost surely invertible. Then, 
\begin{equation}
\boldsymbol{\mathsf{MMSE}}(\X| \Y)= \E \left[  \left( \boldsymbol{\mathsf{J}}_\mathbf{\Y}  {\bf T}(\Y)  \right)^{-1}     \boldsymbol{\mathsf{J}}_\mathbf{\Y} \E[\X|\Y]  \right ]. 
\end{equation} 
\end{corollary} 
\begin{IEEEproof}
The proof follows by observing that  $\boldsymbol{\mathsf{MMSE}}(\X| \Y)=\E[\boldsymbol{\mathsf{Var}}(  \X  | \Y)  ]$ and using~\eqref{eq:Variance_Identity}.
\end{IEEEproof}

\subsection{Recursive Identities}
\label{sec:Recursive_Identites} 
We here seek to express the $(\ell+1)$-th conditional expectation, with $\ell \geq 1$, as a function of the $\ell=1$-st conditional expectation. 
Towards this end, we focus on the univariate case, and we set $U=X^\ell$ in Theorem~\ref{thm:MainIde_Exponential}. 
By doing this, we arrive at the following recursive identity.
\begin{prop}
\label{prop:Recursion}
Assume that $f_{Y|X}(\cdot| \cdot)$ is as in~\eqref{eq:ExpFam}, and let $\E [X^{\ell}|Y=y] = F_{\ell}(y)$.
Then, the following recursive expression holds for all $\ell \geq 1$,
\begin{align}
F_{\ell+1}(y) = \frac{1}{T'(y)} F'_{\ell}(y) +F_1(y) F_{\ell}(y). \label{eq:Recursion_moments}
\end{align}
\end{prop}

A version of the  recursive identity in~\eqref{eq:Recursion_moments} has appeared in the past in the context of a Gaussian noise channel~\cite{jaffer1972note}. In~\cite{DytsoISIT2020extended}, always in the context of a Gaussian noise channel, the authors showed that the identity in~\eqref{eq:Recursion_moments} has several equivalent representations.  

Next, we show an equivalent version of the identity in  Proposition~\ref{prop:Recursion} by solving the recursion in \eqref{eq:Recursion_moments}. The new equivalent version  establishes an expression for the $(\ell+1)$-th conditional expectation, with $\ell \geq 1$, as a function of the $\ell=1$-st conditional expectation. 
\begin{theorem}
\label{thm:RecSolved}
Assume that $f_{Y|X}(\cdot| \cdot)$ is as in~\eqref{eq:ExpFam}, and
define the following operator for $\ell \in \mathbb{N}$,\footnote{Note that, since we assume that $T(y)$ is an analytic function, the derivatives of $T(y)$ have only isolated zeros. Therefore, the set of $y$'s for which  the operator $D^{(\ell)}_y$ is not well-defined is at most a set of measure zero.   }
\begin{equation}
\label{eq:Doperator}
D^{(\ell)}_y =   \underbrace{ \frac{1}{ T'(y)}  \frac{{\rm d}}{ {\rm d} y}  \frac{1}{ T'(y)}  \frac{{\rm d}}{ {\rm d} y}  \dots   \frac{1}{ T'(y)}  \frac{{\rm d}}{ {\rm d} y}}_{ \text{ $\frac{1}{ T'(y)}  \frac{{\rm d}}{ {\rm d} y}$  $\ell$ times}},
\end{equation}
where\footnote{Note that, when the sufficient statistics is a linear function (e.g., Gaussian noise), the operator $D_y^{(\ell)}$ becomes an $\ell$-th derivative.} $D^{(0)}_y =1$. Then, for every $a \in \mathcal{Y}$, we have that
\begin{align*}
&\E [X^{\ell+1}|Y=y] =
\\& \quad \eu^{- \int_{a}^{y} T'(u) \E[X|Y=u]  {\rm{d}}u} D^{(\ell+1)}_y \eu^{\int_{a}^{y} T'(u) \E[X|Y=u]  {\rm{d}}u}.
\end{align*}
\end{theorem}

\begin{IEEEproof}
We start by defining
\begin{align}
\label{eq:Intrmediate1}
g_{\ell}(y) = F_{\ell}(y) \eu^{\int_{a}^{y} T'(u)F_1(u) {\rm{d}}u},
\end{align}
from which we get
\begin{align}
\label{eq:Intrmediate2}
g'_{\ell}(y) =& F'_\ell(y) \eu^{\int_{a}^{y} T'(u) F_1(u)  {\rm{d}}u}  \nonumber
\\& + F_{\ell}(y) \eu^{\int_{a}^{y} T'(u) F_1(u)  {\rm{d}}u} T'(y) F_1(y).
\end{align}
Moreover, from the result in Proposition~\ref{prop:Recursion}, we have that
\begin{align}
\label{eq:Intrmediate3}
T'(y) F_{\ell+1}(y) =  F'_{\ell}(y) +T'(y) F_1(y) F_{\ell}(y) .
\end{align}
By multiplying both sides of~\eqref{eq:Intrmediate3} by $\eu^{\int_{a}^{y} T'(u) F_1(u)  {\rm{d}}u}$ and substituting the expressions in~\eqref{eq:Intrmediate1} and~\eqref{eq:Intrmediate2}, we obtain
\begin{align}
\label{eq:Intrmediate4}
g_{\ell+1}(y) T'(y) = g'_{\ell}(y) &\Rightarrow g_{\ell+1} (y) = \frac{g'_{\ell}(y)}{T'(y)} \nonumber
\\& \stackrel{\eqref{eq:Doperator}}{\Rightarrow} g_{\ell+1} (y) = D_y^{(\ell)} g_{1}(y).
\end{align}
Substituting the definition of $g_{\ell}(y)$ in~\eqref{eq:Intrmediate1} inside~\eqref{eq:Intrmediate4}, we arrive at
\begin{align*}
F_{\ell+1}(y) \eu^{\int_{a}^{y} T'(u) F_1(u)  {\rm{d}}u} = D_y^{(\ell)}F_{1}(y) \eu^{\int_{a}^{y} T'(u) F_1(u)  {\rm{d}}u},
\end{align*}
and, since $\E [X^{\ell}|Y=y] = F_{\ell}(y)$, we obtain
\begin{align*}
\E[X^{\ell+1}|Y=y] =&  \eu^{-\int_{a}^{y} T'(u) \E[X|Y=u]  {\rm{d}}u} 
\\& \cdot D_y^{(\ell)} \E[X|Y=y] \eu^{\int_{a}^{y} T'(u) \E[X|Y=u]  {\rm{d}}u}.
\end{align*}
Finally, note that
\begin{align*}
&\E[X|Y=y] \eu^{\int_{a}^{y} T'(u) \E[X|Y=u]  {\rm{d}}u} 
\\& = \frac{1}{T'(y)} \frac{{\rm{d}}}{{\rm{d}}y} \eu^{\int_{a}^{y} T'(u) \E[X|Y=u]  {\rm{d}}u}.
\end{align*}
This concludes the proof of Theorem~\ref{thm:RecSolved}.
\end{IEEEproof}

We now conclude this discussion with a couple of remarks that point out some implications of the result in Theorem~\ref{thm:RecSolved}.
\begin{remark}
The univariate case considered above can be extended to the multivariate case by setting  $\U= (\X \X^\mathsf{T})^\ell \X, \ell \in \mathbb{N}$ in Theorem~\ref{thm:MainIde_Exponential}, where for a square matrix $\mathbf{A}$ and $\ell \in \mathbb{N}$ we have that $\mathbf{A}^\ell = \underbrace{\mathbf{A} \times \mathbf{A} \times \ldots \times \mathbf{A}}_{\ell \ \text{times}}$.
\end{remark}

\begin{remark}
By using the Tweedy's formula in~\eqref{eq:Tweedy_Formula},  for $\ell \in \mathbb{N}$ and $y\in \mathcal{Y}$,  inside the right-hand side of the expression in Theorem~\ref{thm:RecSolved}, we obtain
\begin{align*}
&\E [X^{\ell+1}|Y=y] \notag
\\& = \eu^{- \int_{a}^{y} \frac{{\rm{d}}}{{\rm{d}}u} \log \frac{f_Y(u)}{h(u)}  {\rm{d}}u} D^{(\ell+1)}_y \eu^{\int_{a}^{y}  \frac{{\rm{d}}}{{\rm{d}}u} \log \frac{f_Y(u)}{h(u)}  {\rm{d}}u} \notag
\\& = \eu^{-\log \frac{f_Y(y)}{h(y)}+\log \frac{f_Y(a)}{h(a)}} D^{(\ell+1)}_y \eu^{\log \frac{f_Y(y)}{h(y)}-\log \frac{f_Y(a)}{h(a)}} \notag
\\& = \eu^{-\log \frac{f_Y(y)}{h(y)}} D^{(\ell+1)}_y \eu^{\log \frac{f_Y(y)}{h(y)}} \notag
\\& =  \frac{h(y)}{f_Y(y)} D^{(\ell+1)}_y \frac{f_Y(y)}{h(y)}.
\end{align*}
\end{remark}

%

\subsection{Identity for Cumulants} 
\label{sec:cumulants}
We here establish fundamental connections between conditional cumulants and conditional moments  for the exponential family defined in~\eqref{eq:ExpFam}.
In particular, for ease of explanation, we focus on the univariate case.
Consider the  \emph{conditional cumulant-generating function},
\begin{equation}
\label{eq:CondCmGenFunc}
K_X(t|Y\!=\!y)\!=\!\log\left( \E[\eu^{tX}| Y\!=\!y] \right),    y \! \in \! \mathcal{Y} \! \subseteq \! \mathbb{R},  t \in \mathbb{R}. 
\end{equation}
The {\em $\ell$-th conditional cumulant} is given by 
\begin{equation}
\label{eq:CondCum}
\kappa_{X|Y=y}(\ell) =\frac{{\rm d}^{\ell}}{ {\rm d} t^{\ell} } K_X(t|Y=y) \Big |_{t=0}, \  \ell \in \mathbb{N}.
\end{equation} 
The next theorem, the proof of which can be found in Appendix~\ref{app:Cumulants}, provides a relationship that we will leverage in Proposition~\ref{prop:Cumulants} to establish a fundamental connection between conditional cumulants and conditional moments.
\begin{theorem}
\label{thm:Cumulants}
Assume that $f_{Y|X}(\cdot| \cdot)$ is as in~\eqref{eq:ExpFam}, and let $D^{(\ell)}_y$ be defined as in~\eqref{eq:Doperator} with $D^{(0)}_y =1$. Then, for $y\in \mathcal{Y}$ and $t \in \mathbb{R}$, we have that
\begin{align*}
 \frac{{\rm d}^\ell}{ {\rm d} t^\ell}  K_{X}( t |Y=y) = &D^{(\ell)}_y  K_{X}(t |Y=y)
 \\&+ D^{(\ell-1)}_y  \E[ X |Y=y].
\end{align*}
\end{theorem} 
%
%
%
\begin{prop} 
\label{prop:Cumulants}
Assume that $f_{Y|X}(\cdot| \cdot)$ is as in~\eqref{eq:ExpFam}.
For $\ell \in \mathbb{N}$ and $y\in \mathcal{Y}$, we have that
\begin{equation}
\kappa_{X|Y=y}(\ell)  =  D^{(\ell-1)}_y  \E[ X |Y=y].  \label{eq:Derivatives_of_CE}
 \end{equation} 
\end{prop} 
\begin{IEEEproof} The proof follows by leveraging the result in Theorem~\ref{thm:Cumulants} and demonstrating that $D^{(\ell)}_y  K_{X}(t |Y=y)|_{t=0}=0$. 
Towards this end, recall the following power series representation of the cumulant-generating function around $t=0$,
\begin{align*}
 &K_{X}(t |Y=y)  \notag\\
 &=  \kappa_{X|Y=y}(1) t+ \frac{\kappa_{X|Y=y}(2) t^2}{2!} + \frac{\kappa_{X|Y=y}(3) t^3}{3!} + \ldots
\end{align*}
Consequently, since the operator $D^{(\ell-1)}_y$ acts only with respect to $y$, we have that $D^{(\ell)}_y  K_{X}(t |Y=y)|_{t=0}=0$. This concludes the proof of Proposition~\ref{prop:Cumulants}.
\end{IEEEproof}
The identity in \eqref{eq:Derivatives_of_CE} has been previously shown in the context of Gaussian noise for $y=0$~\cite{efron2005local} and for all $y\in \mathbb{R}$~\cite{DytsoISIT2020extended}. The work in~\cite{DytsoISIT2020extended} also contains vector generalizations of~\eqref{eq:Derivatives_of_CE}. 

The following concluding remark points out to some interesting consequences of Proposition~\ref{prop:Cumulants}.
\begin{remark}
The result in Proposition~\ref{prop:Cumulants} establishes a novel relationship between conditional cumulants for the exponential family, namely for $\ell \in \mathbb{N}$ and $y\in \mathcal{Y}$, we have that
\begin{equation*}
\kappa_{X|Y=y}(\ell+1)= \frac{1}{ T'(y)}  \frac{{\rm d}}{ {\rm d} y} \kappa_{X|Y=y}(\ell).
\end{equation*}
Moreover, by using the Tweedy's formula in~\eqref{eq:Tweedy_Formula},  for $\ell \in \mathbb{N}$ and $y\in \mathcal{Y}$, we obtain
\begin{equation*}
\kappa_{X|Y=y}(\ell)=D^{(\ell)}_y \log \left( \frac{f_Y(y)}{ h(y)} \right).
\end{equation*}
In other words, the conditional cumulants depend on the joint distribution $P_{X,Y}$ only through the marginal 
$P_Y$.

The result in~\eqref{eq:Derivatives_of_CE} can also be used to find a power series representation of the conditional expectation.  The polynomial approximation of the conditional expectation has received some recent attention and the interested reader is referred to~\cite{ITW2020Dytso} and~\cite{alghamdi2021polynomial}. 
\end{remark}

\appendices
\section{Proof of Lemma~\ref{lem:Prop}}
\label{app:ProofLemmaRelations}
The proof of~\eqref{eq:CE_info_density} follows by using the Markov chain $ \U \leftrightarrow \X \leftrightarrow \Y$ and the definition of information density in~\eqref{eq:InfoDens},
\begin{align*}
\E[\U|\Y={\bf y} ]&=  \int \U \ {\rm d} P_{\U|\X}  \ {\rm d} P_{\X|\Y{\bf= y} }\\
&=  \int \U \ {\rm d} P_{\U|\X}  \ \eu^{\iota_{P_{\X,\Y}}( \X ;{\bf y})} \   {\rm d} P_{\X}\\
&=\E \left[\U  \eu^{\iota_{P_{\X,\Y}}( \X ;{\bf y})} \right]. 
\end{align*}
To show~\eqref{eq:Gradient_Info_density}, observe the following sequence of steps,
\begin{align*}
\nabla_{\bf y}\iota_{P_{\X,\Y}}( {\bf x} ;{\bf y})&= \nabla_{\bf y} \log \frac{f_{\Y|\X}({\bf y}| {\bf x} )}{ f_{\Y}({\bf y})}\\
&= \rho_{\Y|\X}({\bf y}|{\bf x})-\rho_{\Y}({\bf y}). 
\end{align*} 
To show~\eqref{eq:Scond_average}, note that
\begin{align*}
\E [\rho_{\Y|\X}(\Y|\X)|\Y={\bf y} ]& \stackrel{\rm{(a)}}{=} \E \left[\frac{ \nabla_{\bf y} f_{\Y|\X}({\bf y}|  \X )}{f_{\Y|\X}({\bf y}| \X )} |\Y={\bf y} \right]\\
&\stackrel{\rm{(b)}}{=} \E \left[\frac{ \nabla_{\bf y} f_{\Y|\X}({\bf y}|  \X )}{f_{\Y|\X}({\bf y}| \X )}  \eu^{\iota_{P_{\X,\Y}}( \X ;{\bf y})} \right]\\
&\stackrel{\rm{(c)}}{=} \frac{ \E \left[ \nabla_{\bf y} f_{\Y|\X}({\bf y}|  \X ) \right]}{  f_{\Y}({\bf y})} \\
&\stackrel{\rm{(d)}}{=} \frac{ \nabla_{\bf y}  \E \left[  f_{\Y|\X}({\bf y}|  \X ) \right]}{  f_{\Y}({\bf y})} \\
&= \frac{ \nabla_{\bf y}   f_{\Y}({\bf y}) }{  f_{\Y}({\bf y})} \\
&=\rho_{\Y}({\bf y}),
\end{align*} 
where the labeled equalities follow from: $\rm{(a)}$ applying~\eqref{eq:ConScore}; $\rm{(b)}$ using~\eqref{eq:CE_info_density}; $\rm{(c)}$ applying~\eqref{eq:InfoDens}; and $\rm{(d)}$ interchanging the gradient and the expectation that is permissible by using  the Leibniz integral rule, which requires verifying that $\E \left[  \| \nabla_{\bf y} f_{\Y|\X}({\bf y}|  \X ) \| \right]<\infty$ and that  $ f_{\Y|\X}({\bf y}|{\bf x})$ is absolutely continuous in ${\bf y} $. This concludes the proof of Lemma~\ref{lem:Prop}.

\section{Proof of $\rm{A1}$-$\rm{A5}$ in Theorem~\ref{thm:MainIde_Exponential}}
\label{app:AssumExpFamily}
Clearly, $\rm{A1}$ and $\rm{A4}$ in Theorem~\ref{thm:MainResult} and Theorem~\ref{thm:MainIde_Exponential} are the same.
For $\rm{A2}$ in Theorem~\ref{thm:MainResult}, we have that
\begin{align*}
&\E \left[ \left \| \U \right \|   \left \| \nabla_{ \bf y} \iota_{P_{\X,\Y}}( \X ; \Y) \right \| | \Y={\bf y}  \right ]
\\& \stackrel{{\rm{(a)}}}{=} \E \left[ \left \| \U \right \|   \left \|  \rho_{\Y|\X}({\bf y}|{\bf x})-\rho_{\Y}({\bf y}) \right \| | \Y={\bf y}  \right ]
\\& \stackrel{{\rm{(b)}}}{=} \E \left[ \left \| \U \right \|   \left \|  \nabla_{\bf y} \log (h({\bf y})) \!+\! \left( \boldsymbol{\mathsf{J}}_\mathbf{y}  {\bf T}({\bf y})  \right)   \bf{X}-\rho_{\Y}({\bf y}) \right \| | \Y\!=\!{\bf y}  \right ]
\\& \stackrel{{\rm{(c)}}}{=} \E \left[ \left \| \U \right \|   \left \|  g(\bf y) + \left( \boldsymbol{\mathsf{J}}_\mathbf{y}  {\bf T}({\bf y})  \right)   \bf{X} \right \| | \Y={\bf y}  \right ]
\\& \stackrel{{\rm{(d)}}}{\leq} \E \left[ \left \| \U \right \|   \left \|  g(\bf y)  \right \| | \Y={\bf y}  \right ] \!+\!
\E \left[ \left \| \U \right \|   \left \|   \left( \boldsymbol{\mathsf{J}}_\mathbf{y}  {\bf T}({\bf y})  \right)   \bf{X} \right \| | \Y\!=\!{\bf y}  \right ]
\\& \stackrel{{\rm{(e)}}}{\leq} \left \|  g(\bf y)  \right \| \E \left[ \left \| \U  \right \| \!   | \!\Y\!=\!{\bf y}  \right ] \!+\!\left \|    \boldsymbol{\mathsf{J}}_\mathbf{y}  {\bf T}({\bf y})   \right \|_{\star} \!
\E \left[ \left \| \U \right \| \!   \| \bf{X}\| | \Y\!=\!{\bf y}  \right ],
\end{align*}
where the labeled (in)equalities follow from: 
$\rm{(a)}$ using~\eqref{eq:Gradient_Info_density};
$\rm{(b)}$ the fact that, for the exponential family in~\eqref{eq:ExpFam}, the conditional score function is given by 
\begin{align}
 \rho_{\Y|\X}(\Y|\X)&= \nabla_{\bf y} \log (h({\bf y})) +  \nabla_{\bf y} \langle {\bf{x}}, {\bf T}({\bf y}) \rangle \notag \\
& = \nabla_{\bf y} \log (h({\bf y})) + \left( \boldsymbol{\mathsf{J}}_\mathbf{y}  {\bf T}({\bf y})  \right)   \bf{x};
 \end{align}
 $\rm{(c)}$ letting $g({\bf y}) = \nabla_{\bf y} \log (h({\bf y})) -\rho_{\Y}({\bf y})$;
 $\rm{(d)}$ applying the triangle inequality; and
 $\rm{(e)}$ denoting with $\left \|   \boldsymbol{\mathsf{J}}_\mathbf{y}  {\bf T}({\bf y})  \right \|_{\star}$ the operator norm of $\boldsymbol{\mathsf{J}}_\mathbf{y}  {\bf T}({\bf y})$.
 Thus, a sufficient condition to have $\rm{A2}$ in Theorem~\ref{thm:MainResult} satisfied, is to have $\rm{A2}$ and $\rm{A3}$ in Theorem~\ref{thm:MainIde_Exponential} satisfied.

Finally, for $\rm{A3}$ in Theorem~\ref{thm:MainResult}, we have that
\begin{align*}
& \E \left[  \| \nabla_{\bf y} f_{\Y|\X}({\bf y}|  \X ) \| \right]
\\& \stackrel{{\rm{(a)}}}{=} \E \left[  \|  \left( \nabla_{\bf y} \log (h({\bf y})) + \left( \boldsymbol{\mathsf{J}}_\mathbf{y}  {\bf T}({\bf y})  \right)   \bf{X} \right ) f_{\Y|\X}({\bf y}  |  \X ) \|\right ]
\\& \stackrel{{\rm{(b)}}}{\leq} \E \left[ \|  \left( \nabla_{\bf y} \log (h({\bf y})) \right) f_{\Y|\X}({\bf y}|  \X )  \| \right ]
\\& \quad + \E \left[ \| \left( \left( \boldsymbol{\mathsf{J}}_\mathbf{y}  {\bf T}({\bf y})  \right)   \bf{X}\right )  f_{\Y|\X}({\bf y}|  \X )   \|\right ]
\\&\stackrel{{\rm{(c)}}}{=}\E \left[ f_{\Y|\X}({\bf y}|  \X ) \right ]  \|   \nabla_{\bf y} \log (h({\bf y}))     \|
\\& \quad + \E \left[ f_{\Y|\X}({\bf y}|  \X ) \|  \left( \boldsymbol{\mathsf{J}}_\mathbf{y}  {\bf T}({\bf y})  \right)   \bf{X}     \|\right ]
\\&\stackrel{{\rm{(d)}}}{\leq}\E \left[ f_{\Y|\X}({\bf y}|  \X ) \right ]  \|   \nabla_{\bf y} \log (h({\bf y}))     \|
\\& \quad + \|  \left( \boldsymbol{\mathsf{J}}_\mathbf{y}  {\bf T}({\bf y})  \right) \|_{\star} \E \left[ f_{\Y|\X}({\bf y}|  \X ) \|  \bf{X} \|    \|\right ]
\\& \stackrel{{\rm{(e)}}}{=} f_{\Y}({\bf y}) \|   \nabla_{\bf y} \log (h({\bf y}))     \| + \E \left [ \| \X \| | \Y = {\bf y} \right ] f_{\Y}({\bf y}),
\end{align*}
where the labeled (in)equalities follow from: 
$\rm{(a)}$ applying the gradient to $f_{\Y|\X}(\cdot| \cdot)$ in~\eqref{eq:ExpFam};
$\rm{(b)}$ the triangle inequality;
$\rm{(c)}$ the fact that the pdf is non-negative;
$\rm{(d)}$ denoting with $\left \|   \boldsymbol{\mathsf{J}}_\mathbf{y}  {\bf T}({\bf y})  \right \|_{\star}$ the operator norm of $\boldsymbol{\mathsf{J}}_\mathbf{y}  {\bf T}({\bf y})$;
and $\rm{(e)}$ applying Bayes' theorem.
Thus, since $f_{\Y}({\bf y}) \|   \nabla_{\bf y} \log (h({\bf y}))     \|$ is almost surely bounded, a sufficient condition to have $\rm{A3}$ in Theorem~\ref{thm:MainResult} satisfied, is to have $\rm{A5}$ in Theorem~\ref{thm:MainIde_Exponential} satisfied.



\section{Proof of Theorem~\ref{thm:Cumulants}}
\label{app:Cumulants}
Let $U= \eu^{ t X}$. Then, the derivative of the cumulant-generating function can be expressed as,
\begin{align*}
& \frac{{\rm d}}{ {\rm d} y}  K_{X}(t |Y=y)  \notag \\
&\stackrel{{\rm{(a)}}}{=}\frac{ 1  }{ \E[\eu^{ t X}| Y=y ] } \frac{{\rm d}}{ {\rm d} y}  \E[\eu^{ t X}| Y=y ]  \\
&\stackrel{{\rm{(b)}}}{=}\frac{ T'(y)  \  \boldsymbol{\mathsf{Cov}} ( X, \eu^{tX} | Y=y) }{ \E[\eu^{ tX}| Y=y ] } \\
&\stackrel{{\rm{(c)}}}{=}\frac{   T'(y)   \left(  \E[  X \eu^{ tX} | Y=y  ]  -\E[ X |Y=y] \E[  \eu^{t X}  | Y=y]  \right) }{ \E[\eu^{ t X}| Y=y ] }  \\
&{\stackrel{{\rm{(d)}}}{=}\frac{   T'(y)   \left(  \frac{{\rm{d}}}{{\rm{d}}t}\E[  \eu^{ tX} | Y=y  ]  -\E[ X |Y=y] \E[  \eu^{t X}  | Y=y]  \right) }{ \E[\eu^{ t X}| Y=y ] }}  \\
&\stackrel{{\rm{(e)}}}{=} T'(y)  \left(\frac{{\rm d}}{ {\rm d} t}   \log \left(  \E[  \eu^{ t X}  | Y=y] \right)   -\E[ X |Y=y]   \right) \\ 
&\stackrel{{\rm{(f)}}}{=} T'(y)  \frac{{\rm d}}{ {\rm d} t}  K_{X}( t |Y=y)   -    T'(y)     \E[ X |Y=y],
\end{align*}
where the labeled equalities follow from:
$\rm{(a)}$ using the expression in~\eqref{eq:CondCmGenFunc};
$\rm{(b)}$ applying Theorem~\ref{thm:MainIde_Exponential} with $T'(y) = \frac{{\rm{d}}}{{\rm{d}}y}T(y)$;
$\rm{(c)}$ using the definition of conditional cross-covariance in~\eqref{eq:CondCovar};
$\rm{(d)}$ the fact that
\begin{align*}
\E[  X \eu^{ tX} | Y=y  ] = \E \left [ \frac{{\rm{d}}}{{\rm{d}}t} \eu^{tX} | Y=y \right ] = \frac{{\rm{d}}}{{\rm{d}}t}\E[  \eu^{ tX} | Y=y  ],
\end{align*}
where the last equality follows by interchanging the derivative and the expectation, which is permissible by using the Leibniz integral rule under the regularity conditions in Theorem~\ref{thm:MainIde_Exponential};
$\rm{(e)}$ the fact that
\begin{align*}
\frac{\frac{{\rm{d}}}{{\rm{d}}t}\E[  \eu^{ tX} | Y=y  ] }{ \E[\eu^{ t X}| Y=y ]} = \frac{{\rm d}}{ {\rm d} t}   \log \left(  \E[  \eu^{ t X}  | Y=y] \right);
\end{align*}
and $\rm{(f)}$ using the expression in~\eqref{eq:CondCmGenFunc}.
Consequently, we obtain
\begin{equation*}
 \frac{{\rm d}}{ {\rm d} t}  K_{X}( t |Y=y)  = \frac{1}{ T'(y)}  \frac{{\rm d}}{ {\rm d} y}  K_{X}(t |Y=y)+ \E[ X |Y=y].
\end{equation*}
The proof of Theorem~\ref{thm:Cumulants} is concluded by applying induction. 

\bibliography{refs}
\bibliographystyle{IEEEtran}

\end{document}